\documentstyle[12pt,amssymb]{article}

\def\addcontentsline#1#2#3{\relax}

\font\teneufm=eufm10
\font\seveneufm=eufm7
\font\fiveeufm=eufm5
\newfam\eufmfam
\textfont\eufmfam=\teneufm
\scriptfont\eufmfam=\seveneufm
\scriptscriptfont\eufmfam=\fiveeufm

\let\goth\frak
\newcommand{\Z}{{\Bbb Z}}

\newcommand{\C}{{\Bbb C}}

\newcommand{\la}{\lambda}
\newcommand{\La}{\Lambda}
\newcommand{\ve}{\varepsilon}

\newcommand{\ps}{{p^*}} 

\newcommand{\bea}{\begin{eqnarray}}
\newcommand{\ena}{\end{eqnarray}}
\newcommand{\be}{\begin{eqnarray*}}
\newcommand{\en}{\end{eqnarray*}}
\newcommand{\eqref}[1]{(\ref{#1})}

\newcommand\res{\mbox{res}}
\newcommand{\id}{\mbox{id}}
\newcommand{\BW}[5]
{\left(\begin{array}{cc}#1 & #2 \cr #3 & #4 \cr\end{array}
\Biggl| #5\right)}
\newcommand\uq{{U_q\bigl(\widehat{{\goth s}{\goth l}}_2\bigr)}}
\newcommand\slth{{\widehat{{\goth s}{\goth l}}_2}}

\newcommand{\qed}{\hfill \fbox{}\medskip}
\newcommand{\proof}{\medskip\noindent{\it Proof.}\quad }

\newtheorem{thm}{Theorem}[section]
\newtheorem{prop}[thm]{Proposition}
\newtheorem{lem}[thm]{Lemma}

\newcommand{\ignore}[1]{}
\date{\today}      

\begin{document}

%
%

\renewcommand{\thefootnote}{\fnsymbol{footnote}}
\font\csc=cmcsc10 scaled\magstep1

{\baselineskip=14pt
 \rightline{
 \vbox{
       \hbox{September 1997}
}}}

\vskip 11mm
\begin{center}
{\large\bf A coset-type construction for the deformed Virasoro algebra}

\vspace{15mm}

{\csc Michio Jimbo}$\,^{1}$
and
{\csc Jun'ichi Shiraishi}$\,{}^2$
\\ 
%
{\baselineskip=15pt
\it\vskip.35in 
\setcounter{footnote}{0}\renewcommand{\thefootnote}{\arabic{footnote}}
\footnote{e-mail address : jimbo@kusm.kyoto-u.ac.jp}
Division of Mathematics, Graduate School of Science,\\
Kyoto University, Kyoto 606, Japan\\
\vskip.1in 
\footnote{e-mail address : shiraish@momo.issp.u-tokyo.ac.jp}
Institute for Solid State Physics, \\
University of Tokyo, Tokyo 106, Japan \\
}
\end{center}

\begin{abstract}
An analog of the minimal unitary series 
representations for the deformed Virasoro algebra is constructed 
using vertex operators of the quantum affine algebra $\uq$. 
A similar construction is proposed for the elliptic algebra 
${\cal A}_{q,p}\bigl(\slth\bigr)$. 
\end{abstract}


\section{Introduction}

The deformed Virasoro algebra (DVA) introduced in \cite{qVir} 
is presented in terms of a generating series $T(z)=\sum_{n\in\Z}T_nz^{-n}$
as follows. 
\begin{eqnarray}
&&
f\left(\frac{z_2}{z_1}\right) T(z_1)T(z_2)
-T(z_2)T(z_1) f\left(\frac{z_1}{z_2}\right) 
\nonumber\\
&&\qquad =
\frac{(x^{r-1}-x^{-r+1})(x^r-x^{-r})}{x-x^{-1}}
\left(\delta\Bigl(x^{-2}\frac{z_2}{z_1}\Bigr)
-\delta\Bigl(x^{2}\frac{z_2}{z_1}\Bigr)\right),
\label{1.3}\\
&&
f(z)=
\frac{1}{1-z}\frac{(x^{2r}z;x^4)_\infty(x^{-2r+2}z;x^4)_\infty}
{(x^{2r+2}z;x^4)_\infty(x^{-2r+4}z;x^4)_\infty}.
\label{1.4}
\end{eqnarray}
Here $\delta(z)=\sum_{n\in\Z}z^n$ and 
\[
(z;p_1,\cdots,p_k)_\infty=\prod_{n_1,\cdots,n_k=0}^\infty
(1-zp_1^{n_1}\cdots p_k^{n_k}).
\] 
In the limit $x\rightarrow 1$, the Virasoro algebra is recovered 
with the central charge 
\be
c=1-\frac{6}{r(r-1)}.
\en
See also \cite{FR96} 
for a formulation from the viewpoint of chiral vertex algebras. 

A bosonic realization of DVA was presented in \cite{qVir}. 
An analog of the minimal unitary series representations was constructed 
in the work \cite{LukPug2} using the BRST complex. 
(To our knowledge, however, a complete proof about its cohomological 
structure is yet unavailable.) 
Apart from these, little is known about representations of DVA. 
In this note we give an alternative construction of the latter  
representations based on vertex operators (VO's) of 
the quantum affine algebra $\uq$. 
The method is elementary as explained below. 

Recall the analogs of the simplest chiral primary fields 
$\phi_{12}(z)$, $\phi_{21}(z)$ introduced in \cite{LukPug2}:
$\phi_{21}(z)$ is the type I VO's (half transfer matrices), 
while $\phi_{12}(z)$ is the type II VO's 
(particle creation/annihilation operators).
(In the body of the 
text, we will employ the notation $\psi(z)$ for $\phi_{12}(z)$.)
In the conformal case, we have the well known fusion rule 
\bea
\phi_{12}\phi_{12}\sim I + \phi_{13}. 
\label{1.5}
\ena
In the deformed case, the singularity of 
$\rho(z_1/z_2)^{-1}\phi_{12}(z_1)\phi_{12}(z_2)$ 
(with an appropriate scalar factor  $\rho(z)$ \eqref{rho})  
consists of a series of poles 
\be
\frac{z_2}{z_1}=x^2, x^{-2r+4},x^{-4r+6},\cdots
\en
As an analog of \eqref{1.5}, we have 
\bea
\frac{1}{\rho\left(z_1/z_2\right)}
\phi_{12}(z_1)\phi_{12}(z_2)
=
\frac{\id}{1-x^{-2}z_2/z_1} +O(1)
\qquad (z_1\rightarrow x^{-2}z_2).
\label{1.2}
\ena
The DVA current $T(z)$ appears at the second pole $z_2/z_1=x^{-2r+4}$
\cite{8VLuk,JKM96,FR96}:
\begin{equation}
\frac{1}{\rho\left(z_1/z_2\right)}
\phi_{12}(z_1)\phi_{12}(z_2)
=const. \frac{T(z)}{1-x^{2r-4}z_2/z_1}+O(1)
\qquad (z_1\rightarrow x^{-2r+4}z_2).
\label{1.1}
\end{equation}
This fact is naturally expected from \eqref{1.2}, since 
in the conformal case 
the descendant of the identity operator is the Virasoro current itself. 

The formula \eqref{1.1} is easy to establish directly 
in the bosonic realization. 
Alternatively, we can regard \eqref{1.1} as defining the DVA current
in terms of $\phi_{12}(z)$. 
As we will show, 
the commutation relations \eqref{1.3} can be extracted from those of 
$\phi_{12}(z)$ and the structure of its poles, 
without invoking the explicit bosonic expression of $T(z)$. 
In \cite{JMOh} was given a construction of $\phi_{12}(z)$ 
which utilizes the $q$-VO's for the 
quantized affine algebra $\uq$. 
The above procedure then gives rise to a representation of DVA
on the tensor product of integrable modules of $\uq$ 
(see Proposition \ref{prop:c}). 
In the simplest case, this affords 
$c=1/2$ representations of DVA in terms of
a quadratic expression of free fermions. 
The present method is similar in spirit to the coset construction 
of the Virasoro algebra. 
It should be mentioned, however, that the analog of 
the Sugawara operators for DVA is unknown. 

The same recipe can be applied to the VO's for the 
elliptic algebra ${\cal A}_{q,p}(\slth)$, 
assuming their existence and expected analyticity properties. 
In particular, at the Ising point,  
this leads to a fermionic representation of 
DVA at $c=-2$. 
We will comment on these in the last section.

\section{Vertex operators for $\uq$}
\subsection{Notation}

Let $\La_0,\La_1$ be the fundamental weights for $\uq$, and set 
$\rho=\La_0+\La_1$. 
A dominant integral weight of level $k$ has the form 
\bea
\la=\la_l=(k+2-l)\La_0+l \La_1-\rho, 
\qquad 1\le l\le k+1.
\label{la}
\ena
The corresponding integrable highest weight module 
is denoted by  $V(\la)$. 
Let further $V_z$ denote the evaluation module based on
a finite dimensional module $V$. 

The VO's of type I and type II are intertwiners of $\uq$ modules
\begin{eqnarray}
&&\tilde{\Phi}^{(\mu,\la)}(z): V(\la)
\longrightarrow V(\mu)\otimes V_z,
\\
&&\tilde{\Psi}^{*(\eta,\xi)}(z)
: V_z\otimes V(\xi) \longrightarrow V(\eta). 
\end{eqnarray}
We set 
\be
&&
\Phi^{(\mu,\la)}(z)
=\tilde{\Phi}^{(\mu,\la)}(z)z^{\Delta_\mu-\Delta_\la},
\qquad 
\Psi^{*(\eta,\xi)}(z)
=\tilde{\Psi}^{*(\eta,\xi)}(z)z^{\Delta_\eta-\Delta_\xi}. 
\en
For a weight \eqref{la}, $\Delta_\la$ is given by 
\[
\Delta_{\la_l}=\frac{(\la_l,\la_l+2\rho)}{2(k+2)}
=\frac{l^2-1}{4(k+2)}.
\]
We will mostly follow the conventions of \cite{JMOh}, wherein 
$\Phi^{(\mu,\la)}(z)$ and $\Psi^{*(\eta,\xi)}(z)$ are written
as $\Phi^{\mu V}_\la(z)$ and $\Phi^\eta_{V\xi}(z)$, respectively. 
Sometimes we drop the superscripts and write 
 $\Phi(z), \Psi^*(z)$, regarding them as 
operators on the direct sum $\oplus_{l=1}^{k+1} V(\la_l)$. 

In what follows we shall focus attention to the case
\be
V=\C^2, \quad \quad \xi,\eta \mbox{ has level one}.
\en
The components of VO's with respect to the natural basis 
$v_+,v_-\in V$ are defined by 
$\Phi(z) 
=\sum_\ve \Phi_\ve(z)\otimes v_\ve$ and 
$\Psi_\ve^{*}(z)
=\Psi^{*}(z)\left(v_\ve \otimes\cdot  \right)$.
We set $Pv\otimes v'=v'\otimes v$, and 
denote by
$\overline{R}(z)$ the standard $R$-matrix for 
the two-dimensional module (e.g. p.70 in \cite{JM}). 

\subsection{Properties of $\Phi(z),\Psi^*(z)$}

We set
\be
&&r=k+3, \qquad x=-q, \qquad \ps=x^{2(r-1)},
\qquad \{z\}=(z;x^4,\ps)_\infty,
\\
&&\eta_I(z)=\frac{\{\ps x^2 z\}^2}{\{\ps x^4 z\}\{\ps z\}}, 
\qquad
\eta_{II}(z)=\frac{(x^2z;x^4)_\infty}{(z;x^4)_\infty},
\\
&&\rho_{I}(z)=z^{\frac{1}{2r-2}}\frac{\eta_I(z)}{\eta_I(z^{-1})},
\qquad
\rho_{II}(z)=z^{\frac{1}{2}}\frac{\eta_{II}(z)}{\eta_{II}(z^{-1})}.
\en

{}From the theory of $q$-KZ equation we know the following facts.
\begin{enumerate}
\item [(a1)] The product 
\bea 
&&\prod_{i<j}\eta_I(z_j/z_i)\cdot \Phi(z_1)\cdots\Phi(z_N)
\label{ana1}
\ena
is meromorphic
\footnote{By abuse of language, 
we say a function is `meromorphic' if it is a product of 
a meromorphic function and 
power functions in the coordinates $z_1,\cdots,z_N$.}, 
the only singularities being simple poles at 
$z_j/z_i=\ps^{-s} x^{2}$, $i<j, s=1,2,\cdots$. 
Likewise 
\bea 
&&\prod_{i<j}\eta_{II}(z_j/z_i)\cdot \Psi^*(z_1)\cdots\Psi^*(z_N)
\label{ana2}
\ena
is meromorphic, the only singularities being simple poles at 
$z_j/z_i=x^{2}$, $i<j$.
\item[(a2)] 
\bea
&&P\overline{R}(z_1/z_2)
\Phi^{(\nu,\mu)}(z_1)\Phi^{(\mu,\la)}(z_2)
\label{com1}
\\
&&=
\rho_I\left(\frac{z_1}{z_2}\right)
\sum_{\mu'}
\Phi^{(\nu,\mu')}(z_2)\Phi^{(\mu',\la)}(z_1)
W\BW{\la}{\mu}{\mu'}{\nu}{\frac{z_1}{z_2}}, 
\nonumber
\\
&&
\Psi^{*}(z_1)\Psi^{*}(z_2)\left(P\overline{R}(z_1/z_2)\right)^{-1}
=
\rho_{II}\left(\frac{z_1}{z_2}\right)
\Psi^{*}(z_2)\Psi^{*}(z_1). 
\label{com2}
\ena
\item[(a3)] 
\bea
&&
\sum_\ve x^{-\ve/2}
\Phi^{(\nu,\mu)}_\ve(x^{-2}z)\Phi^{(\mu,\la)}_{-\ve}(z)
=x^{2(\Delta_\mu-\Delta_\la)-1/2}g^\mu_\la\delta_{\nu\la}\times\id,
\label{res1}
\\
&&
\res_{z_1=x^{-2}z_2}\Psi^{*}_{\ve_1}(z_1)\Psi^{*}_{\ve_2}(z_2)\frac{dz_1}{z_1}
=
gx^{\ve_2/2}\delta_{\ve_1+\ve_2,0}
\times\id,
\label{res2} 
\ena
where $g=(x^2;x^4)_\infty/(x^4;x^4)_\infty$  and 
$g^\mu_\la$ are constants given in eq.(B3), 
\cite{JMOh} (there $r_\mp$ is misprinted as $r_\pm$). 
\end{enumerate}

In \eqref{com1}, 
\bea
W{\BW{\la}{\mu}{\mu'}{\nu}{z}}
=
\overline{W}^1_k{\BW{\la}{\mu}{\mu'}{\nu}{z}}
\times z^{\Delta_\la+\Delta_\nu-\Delta_{\mu'}-\Delta_\mu-1/{2(k+2)}}
\label{BW}
\ena
denotes the Boltzmann weights for the RSOS model
(see AppendixB, eq.(B.2) in \cite{JMOh}). 
We shall be concerned with the following properties 
rather than their explicit formulas.
\begin{enumerate}
\item $W\BW{\la}{\mu}{\mu'}{\nu}{z}$ is meromorphic on $\C\backslash\{0\}$ with 
simple poles at $z=x^{-2}\ps^s$ ($s\in \Z$),
\item We have the periodicity  
\bea 
&&W\BW{\la}{\mu}{\mu'}{\nu}{\ps z}=W\BW{\la}{\mu}{\mu'}{\nu}{z},
\label{2.31}
\ena
\item At $z=x^{-2}$ we have
\bea
&&\res_{z=x^{-2}}W\BW{\la}{\mu}{\mu'}{\nu}{z}\frac{dz}{z}
=\delta_{\nu\la}a_{\mu\la}b_{\mu'\la},
\label{2.4}
\ena
with some constants $a_{\mu\la},b_{\mu\la}$. 
\end{enumerate}

Notice that the $W$ factor for $\Psi^*(z)$ is simply a scalar.
This is a reflection of the fact that the two-dimensional module $V$ is 
`perfect' (in the sense of crystal base theory) 
for level one representations. 

Property (a3) follows from the fact that 
$\res_{z=x^{-2}}P\overline{R}(z)dz/z$ is proportional to the projector 
onto the trivial module 
\be
V_{x^{-2}z}\otimes V_z \rightarrow \C.
\en

\section{Coset-type construction}

\subsection{Vertex operators of SOS type}

Following \cite{JMOh}, let us introduce the operator 
$\psi^{(\mu,\la;1-i,i)}(z)$ by the composition 
\be
V(\la)\otimes V(\La_i)
\longrightarrow V(\mu)\otimes V_z\otimes V(\La_i) 
\longrightarrow V(\mu)\otimes V(\La_{1-i}).
\en
Namely we set 
\bea
\psi^{(\mu,\la;1-i,i)}(z) 
&=&\left(\id\otimes\Psi^{*(\La_{1-i},\La_i)}(z)\right)
\circ\left(\Phi^{(\mu,\la)}(z)\otimes\id\right)
\nonumber\\
&=&\sum_{\ve=\pm}
\Phi^{(\mu,\la)}_\ve(z)\otimes 
\Psi^{*(\La_{1-i},\La_i)}_\ve(z). 
\label{3.3}
\ena
Dropping superscripts we shall 
often write $\psi^{(\mu,\la)}(z)$ or $\psi(z)$. 
Clearly $\psi^{(\mu,\la)}(z)$ commutes with the diagonal action of $\uq$. 

Set
\bea
&&\eta(z)=\eta_I(z)\eta_{II}(z)=
\frac{\{x^2z\}\{\ps x^2z\}}{\{z\}\{\ps x^4 z\}}, 
\label{eta}\\
&&\rho(z)=\rho_I(z)\rho_{II}(z)=
z^{\frac{r}{2r-2}}\frac{\eta(z)}{\eta(z^{-1})}.
\label{rho}
\ena
Then the properties (a1)--(a3) entail the following. 
\begin{prop}\label{prop:a}
\begin{enumerate}
\item[(b1)] The product 
\bea
&&
\prod_{i<j}\eta(z_j/z_i)\cdot \psi(z_1)\cdots\psi(z_N)
\label{2.1}
\ena
is meromorphic, with at most simple poles at $z_j/z_i=\ps^{-s}x^2$ 
($i<j,s\ge 0$). 
\item[(b2)] As meromorphic functions we have 
\bea
&&
\psi^{(\nu,\mu)}(z_1)\psi^{(\mu,\la)}(z_2)
\label{2.3}\\
&&
=\rho\left(\frac{z_1}{z_2}\right)
\sum_{\mu'}
\psi^{(\nu,\mu')}(z_2)\psi^{(\mu',\la)}(z_1)
W\BW{\la}{\mu}{\mu'}{\nu}{\frac{z_1}{z_2}}.
\nonumber
\ena
\item[(b3)] 
\bea
&&\res_{z_1=x^{-2}z_2}
\psi^{(\nu,\mu)}(z_1)\psi^{(\mu,\la)}(z_2)\frac{dz_1}{z_1}
=\delta_{\nu\la}c_{\mu\la}\times \id,
\label{2.2}
\ena
with $c_{\mu\la}=x^{2(\Delta_\mu-\Delta_\la)-1/2}gg^\mu_\la$. 
\end{enumerate}
\end{prop}

Comparing \eqref{2.2} with \eqref{2.3} and \eqref{2.4} 
we find 
\bea
c_{\la_+,\la}:c_{\la_-,\la}=a_{\la_+,\la}:a_{\la_-,\la}
\label{2.6}
\ena
where $\la_\pm=\la\pm(\La_1-\La_0)$.

\subsection{DVA generators}
 
Property (b3) describes the behavior of 
$\psi(z_1)\psi_(z_2)$ at the `first' pole $z_2/z_1=x^2$. 
We now look at the next pole. 
Define
\bea
T^{(\la;i)}_\pm(z)
&=&\res_{z'=z}\frac{1}{\rho\left(\ps x^{-2}\frac{z'}{z}\right)}
\psi^{(\la,\la_\pm;i,1-i)}(x^{r-2}z')
\psi^{(\la\pm,\la;1-i,i)}(x^{-r+2}z)\frac{dz'}{z'}
\nonumber\\
&&\label{2.71}\\
&=&a_{\la_\pm,\la}
\sum_{\mu'}
\psi^{(\la,\mu';i,1-i)}(x^{-r+2}z)\psi^{(\mu',\la;1-i,i)}(x^{r-2}z)
b_{\mu',\la}. 
\label{2.72}
\ena
In the second line we used the periodicity \eqref{2.31} and 
\eqref{2.4}.

{}From eq.\eqref{2.6} we obtain 
\bea
&&T_-^{(\la;i)}(z)=\frac{c_{\la_-,\la}}{c_{\la_+,\la}}T^{(\la;i)}_+(z).
\label{2.8}
\ena

We claim that the properties (b1)-(b3) imply the following. 

\begin{prop}\label{prop:b} 
\begin{enumerate}
\item[(c1)] As meromorphic functions, 
\bea
&&f\left(\frac{z_2}{z_1}\right)T^{(\la;i)}_\pm(z_1)T^{(\la;i)}_\pm(z_2)
=
f\left(\frac{z_1}{z_2}\right)T^{(\la;i)}_\pm(z_2)T^{(\la;i)}_\pm(z_1),
\label{2.9}
\ena
where $f(z)$ denotes the structure function \eqref{1.4} of DVA. 
\item[(c2)]
The left hand side of \eqref{2.9} is holomorphic in the 
neighborhood of $|z_2/z_1|\le x^{-2}$ except for simple poles at 
$z_2/z_1=x^{-2},x^2$. 
\item[(c3)] The residues are given by 
\be
&&
\res_{z_1=x^{-2}z_2}
f\left(\frac{z_2}{z_1}\right)T^{(\la;i)}_\pm(z_1)T^{(\la;i)}_\pm(z_2)
\frac{dz_1}{z_1}
\\
&&\qquad =-
\res_{z_1=x^{2}z_2}f\left(\frac{z_2}{z_1}\right)
T^{(\la;i)}_\pm(z_1)T^{(\la;i)}_\pm(z_2)\frac{dz_1}{z_1}
\\
&&\qquad
=\frac{(x^{r-1}-x^{-r+1})(x^r-x^{-r})}{x-x^{-1}}
(C_\pm^{(\la;i)})^2 \times \id,
\en
where  $C_\pm^{(\la;i)}$ is a constant related to 
the matrix element $\langle T_\pm^{(\la;i)}\rangle$ 
with respect to the highest weight vector 
$|\la\rangle\otimes |\La_i\rangle$ by (see \eqref{<T>} below)
\be
&&\langle T_\pm^{(\la;i)}\rangle= C_\pm^{(\la;i)}\times
(x^{l_i}+x^{-l_i})
\en
with $\la=\la_l$, $l_0=l$ and $l_1=r-1-l$. 
\end{enumerate}
\end{prop}
The proof of Proposition \ref{prop:b} 
will be given in the next section. 

Let us introduce the normalized operator 
\bea
&&T(z)=(C_\pm^{(\la;i)})^{-1} T^{(\la;i)}_\pm(z). 
\label{T}
\ena
Consider the irreducible decomposition 
\be
V(\la)\otimes V(\La_i)=
\oplus_{\nu} V(\nu)\otimes \Omega_{\la,\La_i;\nu}.
\en
In the right hand side, $\nu$ runs over dominant integral weights of 
level $k+1$, and $ \Omega_{\la,\La_i;\nu}$ signifies the space 
of highest weight vectors of that weight. 
Clearly $T(z)$ acts on $ \Omega_{\la,\La_i;\nu}$. 
{}From the assertions (c1)--(c3) we conclude that 

\begin{prop}\label{prop:c} 
The operator $T(z)$ \eqref{2.71}, \eqref{T}
affords a representation of DVA on 
$\Omega_{\la,\La_i;\nu}$.
\end{prop}

\subsection{Example}

As an example, let us take the case $k=1$, $r=4$. 
In this case $\rho(z)$ and the $W$ weight simplify to 
\be
\rho(z)W\BW{\La_i}{\La_{1-i}}{\La_{1-i}}{\La_i}{z}=-1, 
\en
and (b2) becomes 
\be 
\psi(z_1)\psi(z_2)=-\psi(z_2)\psi(z_1).
\en
Moreover the only poles of both sides are $z_2/z_1=x^{-2},x^2$,
with residues proportional to the identity. 
Writing 
\bea
x^i\psi^{(\La_{1-i},\La_i;1-j,j)}(z)
=\cases{
\psi^{NS}(z) & \hbox{$i\equiv j ~\mbox{mod}~{2}$} \cr
\psi^{R}(z)   &  \hbox{$i\not\equiv j ~\mbox{mod}~2$} \cr
}
\label{fermion}
\ena
we see that the Fourier expansion takes the form 
\be 
&&\psi^{NS}(z)=\sum_{n\in\Z+1/2}\psi_n z^{-n},
\qquad \psi^{R}(z)=\sum_{n\in\Z}\psi_n z^{-n}.
\en
In both cases $\psi(z)=\psi^{NS}(z),\psi^{R}(z)$ 
is a free fermion, whose Fourier components satisfy 
\be
\psi_m\psi_n+\psi_n\psi_m=\delta_{m+n,0}(x^{2m}+x^{-2m}).
\en

The DVA current reads
\be
T(z)=\frac{(1-x^6)}{x^2(1+x^2)}\psi(x^{-2}z)\psi(x^2z).
\en
This is a deformation of the well known representation of
the Virasoro algebra at $c=1/2$
in terms of free fermions. 

\section{Proof of Proposition \ref{prop:b}}

\subsection{Proof of (c1)}

Set 
\be
&&F^{(\la_0,\la_1,\la_2,\la_3,\la_4)}(z_1,z_2,z_3,z_4)
\nonumber\\
&&\qquad =\psi^{(\la_0,\la_1)}(z_1)
\psi^{(\la_1,\la_2)}(z_2)
\psi^{(\la_2,\la_3)}(z_3)
\psi^{(\la_3,\la_4)}(z_4).
\label{2.11}
\en
We have the commutation relation
\be
&&F^{(\la_0,\la_1,\la_2,\la_3,\la_4)}(z_1,z_2,z_3,z_4)
\\
&&=
\prod_{i<j}\rho(z_i/z_j)
\sum_{\mu_1,\mu_2,\mu_3}
F^{(\la_0,\mu_1,\mu_2,\mu_3,\la_4)}(z_3,z_4,z_1,z_2)
\\
&&\qquad \times
W\left(
\begin{array}{ccccc}
\la_0&\la_1&\la_2&\la_3&\la_4\\
\la_0&\mu_1&\mu_2&\mu_3&\la_4\\
\end{array}
\Biggl|{z_1\over z_3},{z_2\over z_3},{z_1\over z_4},{z_2\over z_4}\right),
\en
where the last factor means (see the figure below)
\be
&&
\sum_{\mu}
W\BW{\la_4}{\la_3}{\mu_3}{\mu}{\frac{z_2}{z_4}}
W\BW{\la_3}{\la_2}{\mu}{\la_1}{\frac{z_2}{z_3}}
W\BW{\mu_3}{\mu}{\mu_2}{\mu_1}{\frac{z_1}{z_4}}
W\BW{\mu}{\la_1}{\mu_1}{\la_0}{\frac{z_1}{z_3}}.
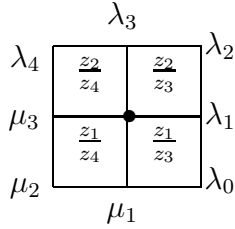
\en
\begin{center}
\be
\mathop{
\raisebox{0pt}[25pt][25pt]
{\mbox{$\mu_3$}}^
{\!\!\!\!\!\!\!\!\mbox{$\la_4$}}_{\!\!\!\!\!\!\!\!\mbox{$\mu_2$}}
{\;\;\;\;\;\;\;\>\,\,\bullet}\!\!\!\!\!\!\!\!\!\!\!\!\!\!\!\!
\begin{array}{|c|c|} \hline
\makebox(18,21){${\frac{z_{2}}{z_{4}}}$} &
\makebox(18,21){${\frac{z_{2}}{z_{3}}}$} \\ \hline
\makebox(18,21){${\frac{z_{1}}{z_{4}}}$} &
\makebox(18,21){${\frac{z_{1}}{z_{3}}}$} \\ \hline
\end{array}
^{\,\mbox{$\la_2$}}_{\,\mbox{$\la_0$}}{\!\!\!\!\!\!\mbox{$\la_1$}}}
^{\mbox{$\la_3$}}_{\mbox{$\mu_1$}}
\en
\end{center}
\medskip

Choosing $(\la_0,\cdots,\la_4)=(\la,\la_\pm,\la,\la_\pm,\la)$, 
taking the residue at $z_1=x^{-2}z_2, z_3=x^{-2}z_4$ and writing 
$z=z_2/z_4$, we find
\bea
&&c_{\la_\pm,\la}^2\times\id
=\rho(z)^2\rho(x^2z)\rho(x^{-2}z)
\nonumber\\
&&\qquad \times
\sum_{\mu,\nu}c_{\mu,\la}c_{\nu,\la}
W\left(
\begin{array}{ccccc}
\la&\la_\pm&\la&\la_\pm&\la\\
\la&\mu&\la&\nu&\la\\
\end{array}
\Biggl|z,x^2z,x^{-2}z,z\right).
\label{2.10}
\ena
Next we take the residue at $z_1=\ps x^{-2}z_2, z_3=\ps x^{-2}z_4$. 
Renaming the variables and using the periodicity \eqref{2.31}
of the Boltzmann weights, we obtain 
\be
&&
T^{(\la)}_\pm(z_1)T^{(\la)}_\pm(z_2)
=\rho(z)^2\rho(\ps^{-1}x^2z)\rho(\ps x^{-2}z)
\\
&&\times
\sum_{\ve_1,\ve_2}
T_{\ve_2}^{(\la)}(z_2)T_{\ve_1}^{(\la)}(z_1)
W\left(
\begin{array}{ccccc}
\la&\la_\pm&\la&\la_\pm&\la\\
\la&\la_{\ve_2}&\la&\la_{\ve_1}&\la\\
\end{array}
\Biggl|z,x^2z,x^{-2}z,z\right)
\en
with $z=z_1/z_2$. 
(Here and after we suppress the index $i$ for brevity.) 
By \eqref{2.8} and \eqref{2.10}
the right hand side simplifies to 
\be
&&
\frac{\rho(\ps^{-1}x^2z)\rho(\ps x^{-2}z)}{\rho(x^2z)\rho(x^{-2}z)}
T^{(\la)}_\pm(z_2)T^{(\la)}_\pm(z_1)
=\frac{f(z_1/z_2)}{f(z_2/z_1)}
T^{(\la)}_\pm(z_2)T^{(\la)}_\pm(z_1),
\en
where we have used the relation 
\bea
&&\frac{\eta(\ps x^{-2}z)\eta(\ps^{-1}x^{2}z)}
{\eta(x^{-2}z)\eta(x^{2}z)}
=
\frac{(1-x^{-2}z)(1-x^2z)}
{(1-\ps^{-1}z)(1-\ps z)}
f(z).
\label{2.13}
\ena
The proof is over. 

\subsection{Proof of (c2)}

We shall show that, in the neighborhood of $|z_2/z_1|\le 1$, 
(i) $f(z_2/z_1)T^{(\la)}_\pm(z_1)T^{(\la)}_\pm(z_2)$ 
is meromorphic, 
(ii) the only pole is $z_2/z_1= x^2$, 
(iii) the pole is simple with residue proportional to the identity. 
The rest of the assertion is then a consequence of the 
commutation relation \eqref{2.9}. 

Recall that 
\be
G^{(\la_0,\la_1,\la_2,\la_3,\la_4)}(z_1,z_2,z_3,z_4)
=\prod_{i<j}\eta(z_j/z_i)
F^{(\la_0,\la_1,\la_2,\la_3,\la_4)}(z_1,z_2,z_3,z_4)
\en 
has poles only at $z_j/z_i=x^2\ps^{-s}$ ($i<j$, $s=0,1,2,\cdots$)
and that they are all simple.  
{}From the definition \eqref{2.72} we see that 
 $T^{(\la)}_\pm(z_1)T^{(\la)}_\pm(z_2)$ 
is a linear combination of the expression
\be
&&\eta\left(\frac{z_2}{z_1}\right)^{-2}
\eta\left(\ps x^{-2}\frac{z_2}{z_1}\right)^{-1}
\eta\left(\ps^{-1}x^{2}\frac{z_2}{z_1}\right)^{-1}
\\
&&\qquad\times
G^{(\la_0,\la_1,\la_2,\la_3,\la_4)}
(x^{-r+2}z_1,x^{r-2}z_1,x^{-r+2}z_2,x^{r-2}z_2).
\en
In the neighborhood of $|z_2/z_1|\le 1$, the last factor can have 
poles only at $z_2/z_1=x^2,\ps,1$, whose multiplicities are at most 
$2,1,1$ respectively.
Multiplying $f(z_2/z_1)$, using \eqref{2.13} and 
\be
&&\eta(z)\eta(x^2z)=\frac{(\ps x^2z;\ps)_\infty}{(z;\ps)_\infty}, 
\en
we find (setting $z=z_2/z_1$) that 
\be
&&f\left(\frac{z_2}{z_1}\right)T^{(\la)}_\pm(z_1)T^{(\la)}_\pm(z_2) 
=
\frac{(\ps x^{-2}z;\ps)_\infty}{(x^2z;\ps)_\infty}
(1-\ps^{-1}z)(1-z)(1-\ps z)
\\
&&\qquad\times 
\hbox{ sum of terms }
G^{(\mu_0,\mu_1,\mu_2,\mu_3,\mu_4)}
(x^{-r+2}z_1,x^{r-2}z_1,x^{-r+2}z_2,x^{r-2}z_2), 
\en
which is regular at $z=\ps,1$. 
The pole $z=x^2$ is apparently double, but the leading term is proportional to 
\be
\frac{1}{(1-x^{-2}z)^2}
\sum_{\mu,\nu}c_{\mu,\la}c_{\nu,\la}
W\BW{\nu}{\la}{\la}{\mu}{x^{-4}}.
\en
Taking the residue of \eqref{2.10} at $z=x^{-2}$ we find that 
the sum is actually $0$.
This shows that the pole is simple. 

\subsection{Proof of (c3)}

The residue can be calculated in the following manner.
\be
&&\res_{z_1=x^{-2}z_2}
f\left(\frac{z_2}{z_1}\right)T^{(\la)}_\pm(z_1)T^{(\la)}_\pm(z_2) 
\frac{dz_1}{z_1}
\\
&&
=\res_{z_1=x^{-2}z_2\atop z_1'=z_1,z_2'=z_2}
\frac{f\left(\frac{z_2}{z_1}\right)}
{\rho\left(\ps x^{-2}\frac{z_1'}{z_1}\right)
\rho\left(\ps x^{-2}\frac{z_2'}{z_2}\right)}
\\
&&\qquad\times
F^{(\la,\la_\pm,\la,\la_\pm,\la)}(x^{r-2}z_1',x^{-r+2}z_1,
x^{r-2}z_2',x^{-r+2}z_2)
\frac{dz_1'}{z_1'}\wedge
\frac{dz_2'}{z_2'}\wedge
\frac{dz_1}{z_1}
\\
&&
=\res_{z_1=x^{-2}z_2\atop z_2'=z_2}
\res_{z_1'=x^{-2}z_2'}
\frac{f\left(\frac{z_2}{z_1}\right)
\rho\left(\ps^{-1}x^{2}\frac{z_1}{z_2'}\right)}
{\rho\left(\ps x^{-2}\frac{z_1'}{z_1}\right)
\rho\left(\ps x^{-2}\frac{z_2'}{z_2}\right)}
\\
&&\qquad\times
\sum_{\mu}
F^{(\la,\la_\pm,\mu,\la_\pm,\la)}(x^{r-2}z_1',x^{r-2}z_2',
x^{-r+2}z_1,x^{-r+2}z_2)
\\
&&\qquad\times
W\BW{\la_\pm}{\la}{\mu}{\la_\pm}{x^2\frac{z_1}{z_2'}}
\frac{dz_1'}{z_1'}\wedge
\frac{dz_2'}{z_2'}\wedge
\frac{dz_1}{z_1}
\en
Since \eqref{BW} is $\delta_{\mu,\mu'}$ at $z=1$, the right hand side becomes 
\be
&&
\res_{z_2'=z_2}
\res_{z_1=x^{-2}z_2'}
\frac{f\left(\frac{z_2}{z_1}\right)
\rho\left(\ps^{-1}x^{2}\frac{z_1}{z_2'}\right)}
{\rho\left(\ps x^{-4}\frac{z_2'}{z_1}\right)
\rho\left(\ps x^{-2}\frac{z_2'}{z_2}\right)}
\\
&&\qquad \times
c_{\la_\pm,\la}\psi^{(\la,\la_\pm)}(x^{-r+2}z_1)
\psi^{(\la_\pm,\la)}(x^{-r+2}z_2)
\frac{dz_2'}{z_2'}\wedge
\frac{dz_1}{z_1}
\\
&&
=-AB c_{\la_\pm,\la}^2
\en
where we have set
\be
&&
A=\res_{z=1}\frac{\rho(\ps^{-1}z)}{\rho(\ps x^{-2}z^{-1})}\frac{dz}{z},
\qquad
B=\frac{f(x^2z)}{\rho(\ps x^{-2}z^{-1})}\Biggl|_{z=1}.
\en

In order to compare this with the matrix element 
 $\langle T_\pm^{(\la;i)}(z) \rangle$, we need the knowledge about 
the two point functions of $\Phi(z),\Psi^*(z)$. 
A formula for the latter can be found e.g. in \cite{JM}.
For $\Phi(z)$, we use the $q$-KZ equation 
\be
&&\langle \Phi^{(\la,\mu)}(\ps z_1)\Phi^{(\mu,\la)}(z_2)\rangle 
\\      
&&\qquad=
x^{\frac{1}{2}}\frac{\eta_I(z_2/z_1)}{\eta_I(\ps^{-1}z_2/z_1)}
(\ps^{-\phi}\otimes 1)\overline{R}(z_1/z_2)
\langle \Phi^{(\la,\mu)}(z_1)\Phi^{(\mu,\la)}(z_2)\rangle 
\en
to reduce the computation to \eqref{res1}
(here $\phi=(\overline{\la}+\overline{\rho})/(r-1)$ and 
$\overline{\la}$ denotes the classical part of $\la$). 
Omitting the details we only give the result.
\bea
&&\langle T^{(\la)}_\pm(z)\rangle 
=(x^{l_i}+x^{-l_i})\times
x^{\frac{r}{2}+\frac{1}{r-1}-\frac{3}{2}}
x^{2(\Delta_{\la_\pm}-\Delta_\la)}g^{\la_\pm}_\la
\nonumber\\
&&\qquad
\times
\frac{(\ps x^4;x^4)_\infty}{(\ps x^2;x^4)_\infty}
\frac{(x^2;\ps)_\infty(\ps x^{-2};\ps)_\infty}
{(\ps;\ps)_\infty^2}.
\label{<T>}
\ena
Here $\la=\la_l$ and $l_0=l, l_1=r-1-l$. 

\section{Elliptic algebra ${\cal A}_{q,p}\bigl(\slth\bigr)$}

We briefly comment on a similar construction of 
DVA current from the level-one VO's of the elliptic algebra
${\cal A}_{q,p}\bigl(\slth\bigr)$. 
(For the notation, see \cite{FIJKMY}.
Here we also write the 
parameters of ${\cal A}_{q,p}\bigl(\slth\bigr)$
as $q=-x$, $p=x^{2r}$ and $\ps=x^{2r-2}$.)
Basically, this complies with Lukyanov's observation \cite{8VLuk}
but with some modification as explained below.

The missing knowledge which is necessary for the 
construction of DVA current from VO's is
the analyticity condition of $n$-point functions 
governed by the $q$-KZ type equation for 
${\cal A}_{q,p}\bigl(\slth\bigr)$.
(See \cite{Fron2} as for the recent progress).
Since quite little is known about this $q$-KZ equation
at present, the best we can do is to start from 
a suitable ansatz for the analyticity. 
Once we assume that, we are able to construct a 
DVA current from the type II VO's as 
\bea
&&T(z)\nonumber\\
&=&C_{II}\res_{\zeta'=\zeta}\sum_\ve
\frac{\xi(\ps^{-1}x^2\zeta^2/\zeta'^2;\ps,q)}
{\xi(\ps x^{-2}\zeta'^2/\zeta^2;\ps,q)}
\Psi^*_\ve(-x^{(r-2)/2}\zeta')
\Psi^*_\ve(x^{-(r-2)/2}\zeta)\frac{d\zeta'}{\zeta'}
\nonumber\\
&=&-C_{II}x^{-r+2}{\Theta_\ps(x^{2}) \over (\ps;\ps)_\infty^3}
\sum_\ve
\Psi^*_\ve(x^{-(r-2)/2}\zeta)
\Psi^*_\ve(-x^{(r-2)/2}\zeta),\label{viraqp}
\ena
satisfying the commutation relation
\begin{eqnarray}
&&
f\left({\zeta_2^2}/{\zeta_1^2}\right) T(\zeta_1)T(\zeta_2)
-T(\zeta_2)T(\zeta_1) f\left({\zeta_1^2}/{\zeta_2^2}\right) 
\nonumber\\
&&\qquad =
{c \over 2}
\left(\delta\Bigl(-x^{-1}\frac{\zeta_2}{\zeta_1}\Bigr)-
\delta\Bigl(x^{-1}\frac{\zeta_2}{\zeta_1}\Bigr)
-\delta\Bigl(-x\frac{\zeta_2}{\zeta_1}\Bigr)+
\delta\Bigl(x\frac{\zeta_2}{\zeta_1}\Bigr)
\right), \label{com}
\end{eqnarray}
where the constant $c$ is chosen as before
\begin{eqnarray}
&&c=
\frac{(x^{r-1}-x^{-r+1})(x^r-x^{-r})}{x-x^{-1}},\nonumber
\end{eqnarray}
and $C_{II}$ is some constant. Note that 
the current $T(\zeta)$ given by (\ref{viraqp}) is 
odd in $\zeta$, i.e. $T(-\zeta)=-T(\zeta)$, and
the
RHS of (\ref{com}) is odd in $\zeta_1$ in accordance with the 
parity of $T(\zeta)$.

At the Ising point $p^{1/2}=x^2$, the type-II VO is realized by
the fermion (\ref{fermion}) as
\be
\sum_\ve \Psi^*_\ve(\zeta)=\psi^{\rm R}(\zeta^2)+\psi^{\rm NS}(\zeta^2),
\en
and 
this gives a fermionic realization of DVA with $c=-2$.
\bea
T(\zeta)
 &=& (x-x^{-1}) \psi^{\rm NS}(\zeta^2)\psi^{\rm R}(\zeta^2).
\ena

We also note that the same can be done with type I VO's and 
observe the composition
\bea
&&C_I\res_{\zeta'=\zeta}\sum_\ve
\frac{\xi(px^2 \zeta'^2/\zeta^2;p,q)}{\xi(p^{-1}x^{-2}\zeta^2/\zeta'^2;p,q)}
\Phi_\ve(-x^{(r+1)/2}\zeta')
\Phi_\ve(x^{-(r+1)/2}\zeta)\frac{d\zeta'}{\zeta'}
\nonumber\\
&=&-C_I x^{r+1}{\Theta_p(x^{-2}) \over (p;p)_\infty^3}
\sum_\ve
\Phi_\ve(x^{-(r+1)/2}\zeta)
\Phi_\ve(-x^{(r+1)/2}\zeta).\nonumber
\ena
enjoyes the DVA relation with a suitable constant $C_I$.
\bigskip

\noindent {\it Acknowledgement.} \quad 
We wish to thank 
Jintai Ding, 
Kenji Iohara, 
Hitoshi Konno, 
Harunobu Kubo, 
Tetsuji Miwa, 
Yas-Hiro Quano
and 
Jun Uchiyama 
for discussions. 


\newpage
%
%
%
\newcommand{\lb}[1]{\label{#1}}
\newcommand{\nn}{\nonumber}
\setcounter{equation}{0}
%
%

\renewcommand{\thefootnote}{\fnsymbol{footnote}}
\font\csc=cmcsc10 scaled\magstep1

{\baselineskip=14pt
 \rightline{
 \vbox{
       \hbox{February 1998}
}}}

\vskip 11mm
\begin{center}
{\large\bf 
Errata to ``Remarks on the deformed Virasoro algebra''}

\vspace{15mm}

{\csc Michio Jimbo}$\,^{1}$
and
{\csc Jun'ichi Shiraishi}$\,{}^2$
\\ 
{\baselineskip=15pt
\it\vskip.35in 
\setcounter{footnote}{0}\renewcommand{\thefootnote}{\arabic{footnote}}
\footnote{e-mail address : jimbo@kusm.kyoto-u.ac.jp}
Division of Mathematics, Graduate School of Science,\\
Kyoto University, Kyoto 606, Japan\\
\vskip.1in 
\footnote{e-mail address : shiraish@momo.issp.u-tokyo.ac.jp}
Institute for Solid State Physics, \\
University of Tokyo, Tokyo 106, Japan \\
}
\end{center}

In the paper \cite{JS}, subsection 4.3, 
the proof of (c3) contains a gap concerning the commutability of 
the operations of taking residues. 
We give here a corrected proof. 

Let 
\bea
&&G^{(\la_0,\cdots,\la_n)}(z_1,\cdots,z_n)=\prod_{i<j}\eta(z_j/z_i)\,
\psi^{(\la_0,\la_1)}(z_1)\cdots \psi^{(\la_{n-1},\la_n)}(z_n).
\lb{G}
\ena
The following are direct consequences of Proposition 3.1.
\begin{lem}\lb{lem1}
\begin{enumerate}
\item[(1)] $G^{(\la_0,\cdots,\la_n)}(z_1,\cdots,z_n)$ 
is holomorphic except for simple poles 
at $z_j/z_i=\ps^{-s}x^2$ ($i<j$, $s\ge 0$).
\item[(2)] We have
\be
&&G^{(\cdots,\la_{k-1},\la_k,\la_{k+1},\cdots)}(\cdots,z_k,z_{k+1},\cdots)
\\
&&\quad
=\left(\frac{z_k}{z_{k+1}}\right)^{\frac{r}{2(r-1)}}
\sum_{\mu}
G^{(\cdots,\la_{k-1},\mu,\la_{k+1},\cdots)}(\cdots,z_{k+1},z_k,\cdots)
W\BW{\la_{k+1}}{\la_k}{\mu}{\la_{k-1}}{\frac{z_k}{z_{k+1}}}.
\en
\item[(3)] As $z_k\rightarrow x^{-2}z_{k+1}$ we have 
\be
&&G^{(\la_0,\cdots,\la_n)}(z_1,\cdots,z_n)=
\frac{\delta_{\la_{k-1}\la_{k+1}}c_{\la_k\la_{k+1}}}{1-x^{-2}z_{k+1}/z_k}
\prod_{i<k}
\frac{(\ps\frac{z_{k+1}}{z_i};\ps)_\infty}
{(x^{-2}\frac{z_{k+1}}{z_i};\ps)_\infty}
\prod_{j>k+1}
\frac{(\ps x^2\frac{z_j}{z_{k+1}};\ps)_\infty}
{(\frac{z_j}{z_{k+1}};\ps)_\infty}
\\
&&\quad \times
\eta(x^2)G^{(\la_0,\cdots,\la_{k-1},\la_{k+2},\cdots,\la_n)}
(z_1,\cdots,z_{k-1},z_{k+2},\cdots,z_n)
+O(1).
\en
\end{enumerate}
\end{lem}

Our aim is to calculate
\bea
&&X=\res_{z_1=x^{-2}z_2}
f\left(\frac{z_2}{z_1}\right)T^{(\la)}_\pm(z_1)T^{(\la)}_\pm(z_2) 
\frac{dz_1}{z_1}.
\lb{X}
\ena
In terms of the function \eqref{G}, this can be written as follows.
\begin{lem}\lb{lem2} We have 
\be
&&X=\res_{z_1=x^{-2}z_2}\varphi(z_2/z_1)\frac{dz_1}{z_1}
\left(\res_{z_1'=z_1}H(z_1',z_1,z_2)\frac{dz_1'}{z_1'}\right),
\en
where
\be
&&H(z_1',z_1,z_2)=
\res_{z_2'=z_2}G^{(\la,\la_\pm,\la,\la_\pm,\la)}
(x^{r-2}z_1',x^{-r+2}z_1,x^{r-2}z_2',x^{-r+2}z_2)
\frac{dz_2'}{z_2'},
\\
&&
\varphi(z)=
\frac{x^{-2r+\frac{2r}{r-1}}}{\eta(\ps x^{-2})^2}
\frac{(\ps x^{-2}z;\ps)_\infty}{(x^{2}z;\ps)_\infty}
(1-\ps^{-1}z)(1-z)(1-\ps z).
\en
\end{lem}

\begin{lem}\lb{lem3}
In the neighborhood of $z_1'=z_1=x^{-2}z_2$, 
$H(z_1',z_1,z_2)$ is holomorphic except for simple poles at 
$z_1'=z_1$, $z_1=x^{-2}z_2$. As $z_1\rightarrow x^{-2}z_2$, we have
\bea
&&H(z_1',z_1,z_2)=-x^{r-2}\frac{c_{\la_\pm,\la}\eta(x^2)}{1-x^{-2}z_2/z_1}
G^{(\la,\la_\pm,\la)}(x^{r-2}z_1',x^{r-2}z_2)
\nn\\
&&\quad\times
\frac{1}{1-\ps}
\frac{(x^2\frac{z_2}{z_1'};\ps)_\infty}
{(\ps^{-1}\frac{z_2}{z_1'};\ps)_\infty}
\frac{(x^2;\ps)_\infty}{(\ps;\ps)_\infty}.
\lb{aa}
\ena
\end{lem}

\proof Consider $G=G^{(\la,\la_\pm,\la,\la_\pm,\la)}
(x^{r-2}z_1',x^{-r+2}z_1,x^{r-2}z_2',x^{-r+2}z_2)$. 
{}From Lemma \ref{lem1}, the only poles of $G$ in the neighborhood of 
$z_1'=z_1=x^{-2}z_2$, $z_2'=z_2$ are
\be
z_1'=z_1, \quad z_2'=z_2,\quad z_1=x^{-2}z_2,\quad z_1'=x^{-2}z_2'.
\en
As $z_1\rightarrow x^{-2}z_2$ we have 
\bea
&&G=\left(\ps^{-1}\frac{z_2}{z_2'}\right)^{\frac{r}{2(r-1)}}
\frac{c_{\la_\pm,\la}}{1-x^{-2}z_2/z_1}
G^{(\la,\la_\pm,\la)}(x^{r-2}z_1',x^{r-2}z_2')
W\BW{\la_\pm}{\la}{\la}{\la_\pm}{\frac{z_2}{z_2'}}
\nn\\
&&\quad\times \eta(x^2)
\frac{(x^2\frac{z_2}{z_1'};\ps)_\infty}
{(\ps^{-1}\frac{z_2}{z_1'};\ps)_\infty}
\frac{(x^2\frac{z_2}{z_2'};\ps)_\infty}
{(\ps^{-1}\frac{z_2}{z_2'};\ps)_\infty}
+O(1).
\lb{bb}
\ena
Likewise, as $z_1'\rightarrow x^{-2}z_2'$ we have 
\bea
&&G=\left(\ps^{-1}x^2\frac{z_1}{z_2'}\right)^{\frac{r}{2(r-1)}}
\frac{c_{\la_\pm,\la}}{1-x^{-2}z_2'/z_1'}
G^{(\la,\la_\pm,\la)}(x^{-r+2}z_1,x^{-r+2}z_2)
W\BW{\la_\pm}{\la}{\la}{\la_\pm}{x^2\frac{z_1}{z_2'}}
\nn\\
&&\quad\times \eta(x^2)
\frac{(x^4\frac{z_1}{z_2'};\ps)_\infty}
{(\ps^{-1}x^2\frac{z_1}{z_2'};\ps)_\infty}
\frac{(x^4\frac{z_2}{z_2'};\ps)_\infty}
{(\ps^{-1}x^2\frac{z_2}{z_2'};\ps)_\infty}+O(1).
\lb{cc}
\ena
Taking the residue of \eqref{cc} at $z_2'=z_2$, we find that 
$H(z_1',z_1,z_2)$ is regular at $z_1'=x^{-2}z_2$. 
The behavior \eqref{aa} is a consequence of \eqref{bb}.
\qed

\medskip

\noindent{\it Proof of (c3).}\quad 
{}From Lemma \ref{lem2} and Lemma \ref{lem3}, 
we can change the order of the residues as
\be
X=\res_{z_1'=x^{-2}z_2}
\frac{dz_1'}{z_1'}
\left(\res_{z_1=x^{-2}z_2}\varphi(z_2/z_1)
H(z_1',z_1,z_2)\frac{dz_1}{z_1}\right).
\en
Using \eqref{aa} we find 
\be
X&=&-x^{r-2}\varphi(x^2)
\res_{z_1'=x^{-2}z_2}\frac{dz_1'}{z_1'}
c_{\la_\pm,\la}\eta(x^2)
\\
&\times&G^{(\la,\la_\pm,\la)}(x^{r-2}z_1',x^{r-2}z_2)
\frac{(x^2\frac{z_2}{z_1'};\ps)_\infty}
{(\ps^{-1}\frac{z_2}{z_1'};\ps)_\infty}
\frac{(x^2;\ps)_\infty}{(\ps;\ps)_\infty}
\\
&=&
-x^{-r+\frac{2}{r-1}}c_{\la_\pm,\la}^2
\left(\frac{\eta(x^2)}{\eta(\ps x^{-2})}\right)^2
(1-x^2)(1-\ps x^2).
\en
After simplification, we obtain the result stated in subsection 4.3.
\qed

\bigskip 

There are also the following corrections. 
\begin{enumerate}
\item The right hand side of Eq.(13) should read
\be
x^{2(\Delta_\mu-\Delta_\lambda)\mp 1/2}
g^\mu_\lambda\delta_{\nu\lambda}\times\id
\qquad \mbox{ for }\mu=\la_\pm. 
\en
\item Line next to Eq.(23), the formula for $c_{\mu\la}$ should read
\be
c_{\mu\la}=x^{2(\Delta_\mu-\Delta_\la)\mp 1/2}gg^\mu_\la
\qquad \mbox{ for }\mu=\la_\pm. 
\en
\item Third line above subsection 4.3: this equation should read
\be
\frac{1}{(1-x^{-2}z)^2}
\sum_{\mu,\nu}b_{\mu,\la}b_{\nu,\la}c_{\mu,\la}c_{\nu,\la}
W\BW{\nu}{\la}{\la}{\mu}{x^{-4}}.
\en
\end{enumerate}

{\it Acknowledgement.}\quad
We thank T. Miwa for critical comments on the paper \cite{JS}.

\end{document}